\documentclass[12pt]{amsart}
\usepackage[T2A]{fontenc}
\usepackage[cp1251]{inputenc}
\usepackage[english]{babel}
\usepackage{amsmath,latexsym,amsthm,amsfonts,tipa,upgreek}
\usepackage{amssymb,amscd,graphpap, stmaryrd}

\newcommand\w{{\omega}}

\newcommand\kk{{\kappa}}
\newcommand\Tau{\mathcal{T}}

\newcommand\FF{{\mathcal F}}

\newtheorem*{Th}{Theorem}

\newtheorem*{Qs}{Question}

\theoremstyle{definition}

\begin{document}

\title{Factoring groups into dense subsets}
\author{Igor Protasov, Serhii Slobodianiuk}
\date{}
\address{Department of Cybernetics, Kyiv University, Volodymyrska 64, 01033, Kyiv, Ukraine}
\email{i.v.protasov@gmail.com; }
\address{Department of Mechanics and Mathematics, Kyiv University, Volodymyrska 64, 01033, Kyiv, Ukraine}
\email{slobodianiuk@yandex.ru}
\maketitle

\begin{abstract} 
Let $G$ be a group of cardinality $\kk>\aleph_0$ endowed with a topology $\Tau$ such that $|U|=\kk$ for every non-empty $U\in\Tau$ and $\Tau$ has a base of cardinality $\kk$. We prove that $G$ could be factorized $G=AB$ (i.e. each $g\in G$ has unique representation $g=ab$, $a\in A$, $b\in B$) into dense subsets $A,B$, $|A|=|B|=\kk$. We do not know if this statement holds for $\kk = \aleph_0$ even if $G$ is a topological group.
\end{abstract}

\section{Introduction}

For a cardinal $\kk$, a topological space $X$ is called {\it $\kk$-resolvable} if $X$ can be partitioned into $\kk$ dense subsets \cite{b1}.  In the case $\kk=2$, these spaces were defined by Hewitt \cite{b4} as {\it resolvable spaces}. If $X$ is not $\kappa$-resolvable then $X$ is called {\it $\kappa$-irresolvable.}

In topological groups, the intensive study of resolvability was initiated by the following remarkable theorem of Comfort and van Mill \cite{b2}:
every countable  non-discrete Abelian topological group $G$ with finite subgroup $B(G)$ of elements  of order 2 is 2-resolvable.
In fact \cite{b11}, every infinite Abelian  group $G$ with finite $B(G)$ can be partitioned into $ \omega$ subsets dense in every non-discrete group topology on $G$. 
On the other hand, under Martin's Axiom, the countable Boolean group $G$, $G=B(G)$ admits maximal (hence, 2-irresolvable)  group topology \cite{b5}.  
Every non-discrete $\omega$-irresolvable topological group $G$ contains an open countable Boolean subgroup provided  that $G$ is Abelian \cite{b6} or countable \cite{b10}, but the existence of non-discrete $\omega$-irresolvable group topology on the countable Boolean group implies that there is a $P$-point in $\omega^{\ast}$ \cite{b6}.  
Thus, in some models of ZFC (see \cite{b8}), every    non-discrete Abelian or countable topological group is $\omega$-resolvable. For systematic exposition of resolvability in topological and left topological group see \cite[Chapter 13]{b3}.

Recently, a new kind resolvability of groups was introduced in \cite{b7}. A group $G$ provided with a topology $\Tau$ is called {\em box $\kk$-resolvable} if there is a factorization $G = AB$ such that $|A|=\kk$ and each subset $aB$ is dense in $\Tau$.
If $G$ is left topological (i.e. each left shift $x\mapsto gx$, $g\in G$ is continuous) then this is equivalent to $B$ is dense in $\Tau$. 
We recall that a product $AB$ of subsets of a group $G$ is {\em factorization} if $G=AB$ and the subsets $\{aB:a\in A\}$ are pairwise disjoint (equivalently, each $g\in G$ has the unique representation $g=ab$, $a\in A$, $b\in B$).
For factorizations of groups into subsets see \cite{b9}. By \cite[Theorem 1]{b7}, if a topological group $G$ contains an injective convergent sequence then $G$ is box $\w$-resolvable.
This note is to find some conditions under which an infinite group $G$ of cardinality $\kk$ provided with the topology could be factorized into two dense subsets of crdinality $\kk$. 
To this goal, we propose a new method of factorization based on filtrations of groups.

\section{Theorem and Question}
We recall that a weight $w(X)$ of a topological space $X$ is the minimal cardinality of bases of the topology $X$.
\begin{Th}\label{t} 
Let $G$ be an infinite group of cardinality $\kk$, $\kk>\aleph_0$, endowed with a topology  $\Tau$ such that $w(G,\Tau)\le\kk$ and $|U|=\kk$ for each non-empty $U\in\Tau$. Then there is a factorization $G=AB$ into dense subsets $A,B$, $|A|=|B|=\kk$.
\end{Th}

We do not know whether Theorem is true for $\kk=\aleph_0$ even if $G$ is a topological group.
\begin{Qs} 
Let $G$ be a non-discrete countable Hausdorff left topological group $G$ of countable weight.
Can $G$ be factorized $G=AB$ into two countable dense subsets?
\end{Qs}
In Comments, we give a positive answer in the following cases: each finitely generated subgroup of $G$ is nowhere dense, the set $\{x^2: x\in U\}$ is infinite for each non-empty open subset of $G$, $G$ is Abelian.

\section{Proof}
We begin with some general constructions of factorizations of a group $G$ via filtrations of $G$.

Let $G$ be a group with the identity $e$ and let $\kk$ be a cardinal.
A family $\{G_\alpha:\alpha<\kk\}$ of subgroups of $G$ is called a {\em filtration} if
\begin{itemize}
\item[(1)] $G_0=\{e\}$, $G=\bigcup_{\alpha<\kk}G_\alpha$;
\item[(2)] $G\alpha\subset G_\beta$ for all $\alpha< \beta$;
\item[(3)] $G_\beta=\bigcup_{\alpha<\beta} G_\alpha$ for every limit ordinal $\beta$.
\end{itemize}
Every ordinal $\alpha<\kk$ has the unique representation $\alpha=\gamma(\alpha)+n(\alpha)$, where $\gamma(\alpha)$ is either limit ordinal or $0$ and $n(\alpha)\in \w$, $\w=\{0,1,\dots\}$.
We partition $\kk$ into two subsets
$$E(\kk) = \{\alpha<\kk:n(\alpha) \text{ is even}\}, \ \ O(\kk) = \{\alpha<\kk: n(\alpha) \text{ is odd}\}.$$

For each $\alpha\in E(\kk)$, we choose some system $L_\alpha$ of representatives of left cosets of $G_{\alpha+1}\setminus G_\alpha$ 
by $G_\alpha$ so $G_{\alpha+1}\setminus G_\alpha = L_\alpha G_\alpha$.
For each $\alpha\in O(\kk)$, we choose some system $R_\alpha$ of representatives of right cosets of $G_{\alpha+1}\setminus G_\alpha$ by $G_\alpha$ so $G_{\alpha+1}\setminus G_\alpha = G_\alpha R_\alpha$.

We take an arbitrary element $g\in G\setminus\{e\}$ and choose the smallest subgroup $G_\gamma$ such that $g\in G_\gamma$.
By $(3)$, $\gamma = \alpha(g)+1$ so $g\in G_{\alpha(g)+1}\setminus G_{\alpha(g)}$.
If $\alpha(g)\in E(\kk)$ then we choose $x_0(g)\in L_\alpha(g)$ and $g_0\in G_\alpha(g)$ such that $g = x_0(g)g_0$.
If $\alpha(g)\in O(\kk)$ then we choose $y_0(g)\in R_\alpha(g)$ and $g_0\in G_\alpha(g)$ such that $g = g_0y_0(g)$.
If $g_0= e$, we make a stop. Otherwise we repeat the argument for $g_0$ and so on.
Since the set of ordinals $<\kk$ is well ordered, after finite number of steps we get the representation
\begin{itemize}
\item[(4)] $$g = x_0(g)x_1(g)\dots x_{\lambda(g)}(g)y_{\rho(g)}\dots y_1(g)y_0(g),$$ 
$$x_i\in L_{\alpha_i(g)},\ \alpha_0(g)>\alpha_1(g)>\dots>\alpha_{\lambda(g)(g)},$$
$$\ y_i\in R_{\beta_i(g)},\ \beta_0(g)>\beta_1(g)>\dots>\beta_{\rho(g)(g)}.$$
\end{itemize}
If either $\{\alpha_0(g),\dots,\alpha_{\lambda(g)}(g)\}=\varnothing$ or $\{\beta_0(g),\dots,\beta_{\rho(g)}(g)\}=\varnothing$ then we write $g = y_{\rho(g)}\dots y_1(g)y_0(g)$ or $g=x_0(g)x_1(g)\dots x_{\lambda(g)}(g)$.
Thus, $G=AB$ where $A$ is the set of all elements of the form $x_0(g)x_1(g)\dots x_{\lambda(g)}$ and $B$ is the set of all elements of the form $y_{\rho(g)}\dots y_1(g)y_0(g)$.
To show that the product $AB$ is a factorization of $G$, we assume that, besides $(4)$, $g$ has a representation 
$$g = z_0z_1\dots z_\lambda t_\rho\dots t_1t_0.$$
If $g\in G_{\alpha+1}\setminus G_\alpha$ and $\alpha\in O(\kappa)$ then $z_0z_1\dots z_\lambda t_\rho\dots t_1\in G_\alpha$ so $t_0=y_0(g)$. If $\alpha\in E(\kappa)$ then $z_1\dots z_\lambda t_\rho\dots t_1t_0\in G_\alpha$ so $z_0=x_0(g)$. We replace $g$ to $gt_0^{-1}$ or to $z_0^{-1}g$ respectively and repeat the same arguments.

Now we are ready to prove Theorem.
Let $\{U_\alpha:\alpha<\kk\}$ be a $\kk$-sequence of non-empty open sets such that each non-empty $U\in \Tau$ contains some $U_\alpha$. Since $|U_\alpha|=\kk$ for every $\alpha<\kk$, we can construct inductively a filtration $\{G_\alpha:\alpha<\kk\}$, $|G_\alpha|=\max\{\aleph_0, |\alpha|\}$ such that, for each $\alpha\in E(\kk)$ (resp. $\alpha\in O(\kk)$) there is a system $L_\alpha$ (resp. $R_\alpha$) of representatives of left (resp. right ) cosets of $G_{\alpha+1}\setminus G_\alpha$ by $G_\alpha$ such that $L_\alpha\cap U_\gamma\neq\varnothing$ (resp. $R_\alpha\cap U_\gamma \neq \varnothing$) for each $\gamma\le\alpha$.
Then the subsets $A,B$ of above factorization of $G$ are dense in $\Tau$ because $L_\alpha\subset A$, $R_\beta\subset B$ for each $\alpha\in E(\kk)$, $\beta\in O(\kk)$.

\section{Comments}
$1.$ Analyzing the proof, we see that Theorem holds under weaker condition: $G$ has a family $\FF$ of subsets such that $|\FF| = \kk$, $|F| = \kk$ for each $F\in\FF$ and, for every non-empty $U\in\Tau$, there is $F\in\FF$ such that  $F\subseteq U$.

If $\kk = \aleph_0$ but each finitely generating subgroup of $G$ is nowhere dense, we can choose a family $\{G_n:n\in\w\}$ such that corresponding $A,B$ are dense.
Thus, we get a positive answer to Question if each finitely generated subgroup $H$ of $G$ is nowhere dense (equivalently the closure of $H$ is not open).

$2.$ Let $G$ be a group and $A, B$ be subsets of $G$. We say that the product $AB$ is a {\em partial factorization} if the subsets $\{aB: a\in A\}$ are pairwise disjoint (equivalently, $\{Ab:b\in B\}$ are pairwise disjoint).

We assume that $AB$ is a partial factorization of $G$ into finite subsets and $X$ be an infinite subset of $G$. Then the following statements are easily verified
\begin{itemize}
\item[(5)] there is $x\in X$ such that $x\notin B$ and $A(B\cup\{x\})$ is a partial factorization;
\item[(6)] if the set $\{x^2:x\in X\}$ is infinite then there is $x\in X$ such that $(A\cup\{x,x^{-1}\})B$ is a partial factorization. 
\end{itemize}

$3.$ Let $G$ be a non-discrete Hausdorff topological group, $AB$ be a partial factorization of 
$G$ into finite subsets, $A=A^{-1}$, $e\in A\cap B$ and $g\notin B$. Then
\begin{itemize}
\item[(7)] there is a neighbourhood $V$ of $e$ such that, for $U=V\setminus\{e\}$ and for any $x\in U$, the product $(A\cup\{x,x^{-1}\})(B\cup\{x^{-1}g\})$ is a partial factorization (so $g\in (A\cup\{x,x^{-1}\})(B\cup\{x^{-1}g\})$).
\end{itemize}
It suffices to choose $V$ so that $V=V^{-1}$ and $AUg\cap AB=\varnothing$, $UB\cap(AB\cup AUg)=\varnothing$, $U^2g\cap AB=\varnothing$, $U\cap A = \varnothing$. 
We use $A=A^{-1}$ only in $UB\cap AUg=\varnothing$.

$4.$ Let $G$ be countable non-discrete Hausdorff topological group such that $\{x^2: x\in U\}$ is infinite for every non-empty open subset $U$ of $G$. We enumerate $G=\{g_n:n\in\w\}$, $g_0 = e$ and choose a countable base $\{U_n:n\in\w\}$ for non-empty open sets. 
We put $A_0=\{e\}$, $B_0=\{e\}$ and use $(5)$, $(6)$, $(7)$ to choose inductively two sequences $(A_n)_{n\in\w}$ and $(B_n)_{n\in\w}$ of finite subsets of $G$ such that for every $n\in\w$
$A_n\subset A_{n+1}$, $B_n\subseteq B_{n+1}$, $A_n=A_n^{-1}$, $A_nB_n$ is a partial factorization,
$g_n\in A_nB_n$, $A_n\cap U_n\neq\varnothing$, $B_n\cap U_n\neq\varnothing$.
We put $A=\bigcup_{n\in\w} A_n$, $B=\bigcup_{n\in\w}B_n$ and note that $AB$ is a factorization of $G$ into dense subsets.

$5.$ Let $G$ be a countable Abelian non-discrete Hausdorff topological group of countable weight.
We suppose that $G$ contains a non-discrete finitely generated subgroup $H$.
Given any non-empty open subset $U$ of $G$, we choose a neighborhood $X$ of $e$ in $H$ and $g\in S$ such that $Xg\subset U$. 
Since $H$ is finitely generated, the set $\{x^2:x\in X\}$ is infinite so we can apply comment $4$. 
If each finitely generated subgroup of $G$ is discrete then, to answer the Question we use comment $1$.

$6.$ Let $G$ be a countable group endowed with a topology $\Tau$ of countable weight such that $U$ is infinite for every $U\in\Tau$. Applying the inductive construction from comment $5$ to $A_nB_n$ and $B_{n+1}^{-1}A_n^{-1}$, we get a partial factorization of $G$ into two dense subsets.

$7.$ Let $G$ be a group satisfying the assumption of Theorem and let $\gamma$ be an infinite cardinal, $\gamma<\kk$. We take a subgroup $A$ of cardinality $\gamma$ and choose inductively a dense set $B$ of representatives of right cosets of $G$ by $A$. Then we get a factorization $G=AB$. In particular, if $G$ is left topological then $G$ is box $\gamma$-resolvable.

\end{document}